\documentclass{elsarticle}

\usepackage{amsmath}
\usepackage{amsfonts} 
\usepackage{bm}
\usepackage{tikz}
\usepackage{pgfplots}
\usepgflibrary{arrows}
\usepackage{rotating} 
\graphicspath{{pictures/}}

\newtheorem{thm}{Theorem}

\newproof{pf}{Proof}

\newcommand{\const}{\mathop{\rm const}\nolimits}

\journal{arXiv} 

\begin{document}

\begin{frontmatter}

\title{Splitting methods for solution decomposition in nonstationary problems}

\author[tamu,uni]{Yalchin Efendiev\corref{cor}}
\ead{efendiev@math.tamu.edu}
\author[nsi,uni]{Petr N. Vabishchevich}
\ead{vabishchevich@gmail.com}

\address[tamu]{Department of Mathematics, Texas A\& M University, College Station, TX 77843, USA}
\address[nsi]{Nuclear Safety Institute, Russian Academy of Sciences, Moscow, Russia}
\address[uni]{North-Eastern Federal University, Yakutsk, Russia}

\cortext[cor]{Corresponding author}

\begin{abstract}
In approximating  solutions of nonstationary problems, various approaches are used to compute the solution at a new time level from a number of simpler (sub-)problems. Among these approaches are splitting methods. Standard splitting schemes are based on one or another additive splitting of the operator into ``simpler'' operators that are more convenient/easier for the computer implementation  and use inhomogeneous (explicitly-implicit) time approximations. In this paper, a new class of splitting schemes is proposed that is characterized by an additive representation of the solution instead of the operator corresponding to the problem (called problem operator).  
A specific feature of the proposed splitting is that 
the resulting coupled equations for individual solution components 
 consist of the time derivatives of the solution components.
The proposed 
approaches are motivated by various applications, including multiscale methods, 
domain decomposition, 
and so on, where spatially local problems are 
solved and used to compute the solution.
 Unconditionally stable splitting schemes are constructed for a first-order evolution equation, which is considered in a finite-dimensional Hilbert space. 
In our splitting algorithms, we consider  the decomposition
 of both the main operator
of the system and the operator at the time derivative. Our goal is to provide
a general framework that combines temporal splitting algorithms and
spatial decomposition and its analysis. Applications of the framework
will be studied separately.

\end{abstract}

\begin{keyword}
first-order evolution equations, Cauchy problem, solution decomposition, 
decomposition methods, splitting methods, stability of the difference
schemes, system of evolution equations.
\MSC[2010] 65J08 \sep 65M06 \sep 65M12
\end{keyword}

\end{frontmatter}

\section{Introduction}

Development of computational algorithms for the approximate solution of 
nonstationary  partial differential equations
is carried out on parabolic and hyperbolic
equations of the second order. 
The issues related to the time approximation require a separate consideration
when using finite element or finite volume approximations in space
\cite{Ascher2008,Gustafsson2008}. 
Explicit schemes \cite{LeVeque2007,SamarskiiTheory}, that are easy to 
implement, have  strict stability constraints for time steps.
Implicit schemes, on the contrary, belong to the class of 
unconditionally stable schemes, but the computation of the solution at new
time step is more complicated compared to explicit schemes. 
The goal of many studies is to build
time approximations that would 
keep advantages of both
 explicit schemes (simple computational implementation),
and implicit schemes (unconditional stability).

Some options for simplifying the computation on a 
new time level (where the solution is computed)
without loss of stability
associated with the use of non-uniform time approximations.
The problem operator is split into two operator terms
by selecting an appropriate decomposition  for 
computational implementation 
purposes
- for example, linear, stationary.
In explicit-implicit schemes (IMEX methods) (see, for example, 
\cite{Ascher1995,HundsdorferVerwer2003}),
one part of the problem operator 
is taken from the lower 
time level, and the second (admissible) part is taken from the upper (current)
level.

The main idea of ​​constructing computationally acceptable unconditionally stable
schemes is implemented when building splitting schemes 
\cite{Marchuk1990,VabishchevichAdditive}.
In this case, the transition to a new level in time is carried out by solving 
evolutionary problems for individual
operator terms. The splitting schemes are characterized by the 
choice of splitting operators and the
formulation of auxiliary problems that are used to
determine an approximate solution.
In the case of two-component splitting, the most promising approaches
 are 
the operator analogs of the classical
ADI (Alternating Direction Implicit) 
schemes \cite{PeacemanRachford1955,DouglasRachford1956}.
Among the unconditionally stable 
multicomponent splitting schemes,
we  note
 the schemes of summarized approximation \cite{SamarskiiTheory,Marchuk1990},
regularized additive schemes \cite{SamarskiiVabischevich1998,Vabischevich2010a}
and vector schemes \cite{Abrashin1990,Vabischevich1994}.

Standard splitting schemes are based on the initial additive 
decomposition of the operator(s) of the problem, that are typically
used for approximation.
An example of such technology can be the additive-averaged schemes 
\cite{GordezianiMeladze1974},
when the solution on the new time level is the arithmetic mean of 
auxiliary problems' solutions \cite{VabishchevichAdditive}.
In many cases, a methodologically more acceptable approach is to consider
the decomposition of the solution, which indirectly
yields a decomposition of the operator.
For example, when constructing domain decomposition 
schemes for the approximate nonstationary solution, 
it is natural to consider isolating the solution in the subdomains 
\cite{ToselliWidlund2005,Mathew2008} and performing 
 solution splitting, as a result.
After that, we select tasks in subdomains, perform operator decompositions, 
and construct certain splitting schemes
\cite{SamarskiiMatusVabischevich2002}.

In the present work, 
the transition to simpler problems is carried out on the basis 
of the solution decomposition. Based on the solution decomposition, 
we construct splitting
schemes for the solution. 
An approximate solution to the Cauchy problem in a finite-dimensional 
Hilbert space for a first-order evolution equation
is constructed on a family of spaces
 using proposed restriction and prolongation operators.
These restrictions and prolongation operators are derived from solution
decomposition.
After that, the individual components of the solution are determined from the system of evolutionary equations.
Three-level schemes for splitting the solution with time transfer to the upper level are proposed and investigated based on splitting that uses
the diagonal parts of the operator matrices of the corresponding evolutionary 
system  for the components of the solution.
When narrowing the class of solution decomposition methods to the case of
direct sum of subspaces, two- and three-level splitting schemes are constructed
that are of the second order of accuracy in time.

Our studies are motivated by a number of problems that include domain decomposition methods \cite{ToselliWidlund2005}, multiscale methods \cite{chung2016adaptive}, reduced-order models, proper orthogonal decomposition techniques \cite{pinnau2008model}, and so on. In these approaches, the solution is decomposed in space using various ideas. For example, in multiscale methods, the decomposition of the solution is based on local basis functions (cf., domain decomposition methods), in proper orthogonal decomposition techniques, the solution decomposition uses global reduced-order basis functions. These decompositions provide appropriate restriction and prolongation operators. When using these techniques for non-stationary problems, one can take an advantage of using implicit time stepping only in some parts of the solution space and, thus, reducing the computational cost at each time iteration. Though these time decompositions may be intuitive, their rigorous analysis and a choice of parameters in time splitting require further investigations. This paper presents a general framework that will be used to guide these spatial decomposition approaches in non-stationary problems. We will report our results regarding applications of the proposed framework elsewhere.

The paper is organized as follows. In Section \ref{sec:2}, we discuss solution splitting and some preliminaries. In Section \ref{sec:3}, we discuss splitting schemes and provide stability conditions. Section \ref{sec:4} is devoted to solution decomposition on direct sum subspaces.

\section{Solution splitting}
\label{sec:2}

Let  $U$ be a finite dimensional Hilbert space. We consider a Cauchy
problem for the first-order evolution equation:
\begin{equation}\label{1}
 \frac{d u}{d t} + A u = f(t),
 \quad 0 < t \leq T, 
\end{equation} 
\begin{equation}\label{2}
 u(0)= u^0 .
\end{equation}
We seek the solution 
$u(t)$ of (\ref{1}) for $0 < t \leq T$ in a finite dimensional
Hilbert space  
$U$ with prescribed right hand side $f(t)$ and the initial condition
(\ref{2}).
For simplicity, we assume  $A: U \mapsto U$ does not depend on $t$, self-adjoint, and positive:
\begin{equation}\label{3}
 \frac{d}{d t} A = A \frac{d}{d t},
 \quad A = A^* >  0 . 
\end{equation} 
The equations (\ref{1})--(\ref{3}) are typically obtained as a result
of space discretization of parabolic partial differential equations. 
When using finite-difference
approximations \cite{SamarskiiTheory}, $u$ is a grid-based 
function defined at the nodes of a 
computational grid. With finite element approximations
\cite{Thomee2006}, the unknowns $u$ are the coefficients 
 in the corresponding finite element basis functions
for the approximate
solutions.

The scalar product for $ u, v \in U $ is $ (u, v) $, and the norm is 
$\| u \| = (u, u) ^ {1/2} $.
For a self-adjoint and positive operator $ D $, the Hilbert space is defined
$ U_D $ with scalar product and norm $ (u, v) _D = (D u, v), \ \| u \| _D = (u, v) _D ^ {1/2} $.

In the approximate solution of the Cauchy problem (\ref{1}), (\ref{2}), 
 implicit time approximations are often used,
which provide unconditionally stable schemes. 
We will use a uniform grid in time with
step $\tau$ and let  $y^n=y(t^n), \ t^n = n\tau$,
$n =0, \dots, N, \ N \tau = T $. 
We can, 
for example, use an implicit scheme:
\begin{equation}\label{4}
 \frac{y^{n+1} - y^{n}}{\tau } + A y^{n+1} = f^{n},
 \quad n = 0,\dots,N-1 , 
\end{equation} 
\begin{equation}\label{5}
 y^0 = u^0 .
\end{equation} 
Difference scheme (\ref{4}), (\ref{5}) approximates (\ref{1}), (\ref{2})
with the first order in $\tau$ assuming 
 sufficient smoothness of the solution $u(t)$.

When formulating stability conditions for two-level schemes, we can focus on
general stability results for operator-difference schemes 
\cite{SamarskiiTheory,SamarskiiMatusVabischevich2002}.
Our goal is to obtain similar  a priori estimates that take 
place for the differential problem.
The simplest estimate of the stability of the solution with respect to 
the initial data and the right-hand side for solving the problem
 (\ref{1}) - (\ref{3})
can be obtained by 
multiplying the equation (\ref{1}) by $du/dt$, which gives
\[
 \left \| \frac{d u}{d t} \right \|^2 +  \frac{1}{2} \frac{d }{d t} (Au,u) = \left (f, \frac{d u}{d t} \right ) .
\] 
Taking into account the positivity of the operator $A$ and 
using the inequality
\[
 \left (f, \frac{d u}{d t} \right ) \leq \left \| \frac{d u}{d t} \right \|^2 + \frac{1}{4} \| f\| ,
\] 
we get
\[
 \frac{d \|u\|_A^2}{d t} \leq \frac{1}{2} \|f\|^2 .
\] 
Applying Gronwall's lemma, we obtain
\begin{equation}\label{6}
 \|u(t)\|_A^2 \leq \|u^0\|_A^2 + \frac{1}{2}\int_{0}^{t} \|f(\theta)\|^2 d \theta .
\end{equation}

For an approximate solution, an analogous process can be done. 
We multiply (\ref{4}) 
 by $y^{n+1} -y^{n}$ and arrive at the equality
\[
 \left \|\frac{y^{n+1}-y^{n}}{\tau }\right \|^2 + (A y^{n+1}, y^{n+1}) = (A y^{n+1},y^{n}) + (f^{n},y^{n+1}-y^{n}). 
\] 
Taking into account
\[
 (A y^{n+1},y^{n}) \leq \frac{1}{2} (A y^{n+1}, y^{n+1}) + \frac{1}{2} (A y^{n}, y^{n}) ,
\]
\[
 (f^{n}, y^{n+1}-y^{n}) \leq \tau \left \|\frac{y^{n+1}-y^{n}}{\tau }\right \|^2+  \frac{1}{4} \tau (f^{n}, f^{n}) ,
\]
we have the following estimate on new time level
\[
 \|y^{n+1}\|_A^2 \leq \|y^{n}\|_A^2 + \frac{1}{2} \tau \|f^{n}\|^2 ,
 \quad n = 0, \dots, N-1 . 
\] 
From here (Gronwall's difference lemma), we obtain an a priori estimate
\begin{equation}\label{7}
 \|y^{n+1}\|_A^2 \leq \|u^0\|_A^2 + \frac{1}{2} \sum_{k = 0}^{n} \tau \|f^{k}\|^2 ,
\end{equation} 
which is a grid analogue of the estimate (\ref{6}).

The main disadvantage of implicit schemes is associated with 
the computational complexity involved in determining
solutions on a new time level.
When using (\ref{4}), (\ref{5}),
 an approximate solution to $y^{n + 1}$ is found as a solution
of the problem
\[
 (I + \tau A) y^{n+1} = \varphi^n ,
\] 
where $ I $ is the unit operator,
for the known right-hand side $\varphi^n = y^{n} + \tau f^{n} $.
A task consists of constructing  such approximations in time, which 
can lead to simpler algorithms
for computing an approximate solution on a new time level.

The splitting methods \cite{Marchuk1990,VabishchevichAdditive} 
use one or another
representation of the operator (operators) of the problem as a sum of simpler operators.
In this case, various variants of explicit-implicit approximations in time are used for
problem (\ref{1}), (\ref{2}), 
when the operator $A$ has an additive representation:
\[
 A = \sum_{i = 1}^{p} A_i .
\] 
The transition to a new time level is associated with the solution of a number of simpler problems
associated with individual operators $ A_i, \ i = 1, \dots, p,$
\[
 (I + \tau \sigma_i A_i) y_i^{n+1} = \varphi_i^n,
 \quad  \sigma_i = \const,
 \quad i = 1, \dots, p .   
\]

We consider a new class of methods that are associated not with 
operator splitting, but with splitting
the solution itself. We 
assume that the solution has an additive representation
\begin{equation}\label{8}
  u = \sum_{i = 1}^{p} u_i .
\end{equation}

The main issue is the construction of individual solution terms  to
achieve the goal of simplifying the tasks for finding these solution terms.
In this paper, it is assumed that this choice problem has been 
solved and we will focus on constructing approximations
in time when solving the Cauchy problem (\ref{1}), (\ref{2}) 
when decomposing the solution in the form (\ref{8})
under the following  general assumptions.

We will work with the family of finite-dimensional Hilbert spaces 
$V_i, \ i = 1, \dots, p $.
For each of these spaces, the linear restriction operator $ R_i $ and 
the interpolation (prolongation or extension) operator $ R^*_i$ are defined:
\[
 R_i: U \mapsto V_i,
 \quad  R^*_i: V_i \mapsto U, 
 \quad i = 1,\dots,p . 
\] 
We assume that for $U$, the following expansion holds
\[
 U = \sum_{i=1}^{p} R^*_i V_i ,
\] 
such that for each 
$u \in U$, the following holds
\begin{equation}\label{9}
 u = \sum_{i=1}^{p} R^*_i v_i ,
 \quad v_i \in V_i,
 \quad i = 1,\dots,p .   
\end{equation} 
In this way, 
$u_i = R^*_i v_i, \ i = 1,\dots,p,$ in (\ref{8}). 

To obtain a system of equations for individual terms 
$v_i, \ i = 1, \dots, p $, we substitute the representation (\ref{9})
into equation (\ref{1}) and multiply it  by $ R_i, \ i = 1, \dots, p $.
This gives
\begin{equation}\label{10}
 \sum_{j=1}^{p} R_i R^*_j \frac{d v_j}{d t}  + \sum_{j=1}^{p} R_i A R^*_j v_j = f_i(t),
 \quad f_i(t) = R_i f(t),
 \quad i = 1,\dots,p .   
\end{equation} 
This system of equations is supplemented by the initial conditions
\begin{equation}\label{11}
 v_i(0) = v_i^0,
 \quad i = 1,\dots,p ,  
\end{equation} 
which follow from (\ref{2}).
In (\ref{11}), $v_i^0, \ i = 1, \dots, p, $ ---
 a solution, possibly not unique, of the system of equations
\[
 \sum_{j=1}^{p} R_i R^*_j v^0_j = R_i u^0,
 \quad i = 1,\dots,p . 
\]

We would like to construct time approximations for the Cauchy problem 
(\ref{10}), (\ref{11}),
which gives the transition to a new level in time by solving 
individual problems for $ v_i ^ 0, \ i = 1, \dots, p $.
In this case, taking into account the representation (\ref{9}), 
we end up with 
the schemes of splitting the solution for the problem (\ref{1}), (\ref{2})

\section{Splitting schemes} 
\label{sec:3}

Let us write the system (\ref{10}) 
in the form of one first-order equation for vector quantities.
Define a vector 
$\bm v = \{v_1, \ldots, v_p \} $ and $\bm f = \{f_1, \ldots, f_p \} $
and from (\ref{10}), (\ref{11}), we get to the Cauchy problem
\begin{equation}\label{12}
 \bm C \frac{d \bm v}{d t} + \bm B \bm v = \bm f ,
\end{equation} 
\begin{equation}\label{13}
 \bm v(t) = \bm v^0 .
\end{equation} 
For the operator matrices $\bm C$ and $\bm B$,
 we have the representation
\[
 \bm C = \{R_i R^*_j \},
 \quad \bm B = \{R_i A R^*_j \} ,
 \quad i, j = 1,\dots,p . 
\]

Problem (\ref{12}), (\ref{13}) can be 
considered on the direct sum of spaces $\bm V = V_1 \oplus \ldots  \oplus V_p$,
when for $\bm v, \bm z \in \bm V$,
 the scalar product and norm are determined by the expressions
\[
 (\bm v, \bm z) = \sum_{1=1}^{p} (v_i, z_i)_i,
 \quad \|\bm v\| =  (\bm v, \bm v)^{1/2} ,
\]
where $(v_i, z_i)_i$ ---scalar product for 
$v_i, z_i \in V_i$ в $V_i, \ i = 1,\dots, p$.

Because
$(R_i R^*_j)^* = R_j R^*_i, \ i, j = 1,\dots,p,$ the operator  $\bm C$ 
is self-adjoint. In addition, we have
\begin{equation}\label{14}
  (\bm C \bm v, \bm v) = \left ( \left (\sum_{i=1}^{p} R^*_i v_i\right )^2, 1 \right ) .
\end{equation}
From (\ref{14}), it follows that the operator
$\bm C$ is non-negative, and taking into account
(\ref{9}), we have
\begin{equation}\label{15}
  (\bm C \bm v, \bm v) = \|u\|^2.
\end{equation}
Similar properties are set for the operator $\bm B$ conditions (\ref{3}):
\begin{equation}\label{16}
  (\bm B \bm v, \bm v) = \left ( \left (\sum_{i=1}^{p} R^*_i v_i\right )^2, 1 \right )_A 
  = \|u\|_A^2.
\end{equation}
We arrive at vector-valued problem
(\ref{12}), (\ref{13}), where
\begin{equation}\label{17}
  \bm C = \bm C^* \ge 0,
  \quad \bm B = \bm B^* \ge 0 .
\end{equation}

Multiplying  
(\ref{12}) in  $\bm V$ by $d \bm v/dt$ (in scalar way), we get
\[
  \left (\bm C \frac{d \bm v}{d t} , \frac{d \bm v}{d t} \right ) + 
  \frac{1}{2} \frac{d }{d t}(\bm B \bm v, \bm v) = \left (\bm f, \frac{d \bm v}{d t} \right)  .
\] 
For the right hand side, we have
\[
 \left (\bm f, \frac{d \bm v}{d t} \right ) = \sum_{i=1}^{p}\left ( R_i f, \frac{d v_i}{d t} \right )_i = 
 \sum_{i=1}^{p} \left (f, R^*_i \frac{d v_i}{d t} \right ) = \left (f,\frac{d u}{d t} \right ) .
\] 
If we take into account
(\ref{15}), (\ref{17}), we get the estimate (\ref{6})
for the solution of 
(\ref{1}), (\ref{2}).

When constructing difference schemes for problem (\ref{12}), (\ref{13}), 
we start from an implicit scheme:
\begin{equation}\label{18}
 \bm C \frac{\bm w^{n+1} - \bm w^{n}}{\tau } + \bm B \bm w^{n+1} = \bm f^{n} ,
 \quad n = 0,\ldots,N-1, 
\end{equation} 
\begin{equation}\label{19}
 \bm w^{0} = \bm v^0 .
\end{equation} 
To obtain an a priori estimate, we  scalarly multiply (\ref{18})  by $\bm w^{n + 1} - \bm w^{n}$.
As in the case of the Cauchy problem for the differential-operator equation 
(\ref{12}), (\ref{13}), for the scheme
(\ref{18}), (\ref{19}), we have the stability estimate (\ref{7}),
which is also the case for the implicit scheme (\ref{4}), (\ref{5}). 
In this case
\[
 y^n = \sum_{i=1}^{p} R^*_i w^n_i ,
 \quad n = 0, \dots,N . 
\]

For the vector $\bm v$, we cannot obtain similar estimates 
in view of the fact that
(see (\ref{17})) neither $\bm C$ nor $\bm B$ are positive definite
and the norm $\bm v$ cannot be associated with them. 
If we consider (\ref{12}), (\ref{13}) as an auxiliary problem
for solving (\ref{1}), (\ref{2}), then we are 
not interested in the vector $\bm v$, but in the scalar 
combination of its components $u$.
The stability estimates, we have
both at the differential and at the discrete level for the computed
 values ​​$u(t)$ and $y^n$.

Implementing an implicit scheme requires solving a vector problem
\begin{equation}\label{20}
 (\bm C + \tau \bm B) \bm w^{n+1} = \bm \varphi^n,
 \quad \bm \varphi^n = \bm w^{n} + \bm f^{n} , 
\end{equation} 
at every time step. In the general case, we have (\ref{17}) and, 
therefore, can not guarantee the uniqueness.
We must narrow the admissible class of
spaces $ V_i, \ i = 1, \dots, p $ (operators $ R_i, \ i = 1, \dots, p $). 
For example,
such that the operator $\bm C$ is positive. 
In this case, instead of (\ref{17}) we have
\begin{equation}\label{21}
  \bm C = \bm C^* > 0,
  \quad \bm B = \bm B^* > 0 .
\end{equation}
The problems of computational implementation lie in the fact that the individual components of the vector on each new time level are determined from
coupled system of equations (\ref{20}). Because of this, we focus on 
the splitting schemes of the solution, which separate components
$v_i, \ i = 1, \dots, p $, of the vector $\bm v$, and  
are determined from independent problems.

We can construct splitting schemes for the solution by separating 
the diagonal parts of the operator matrices $\bm C$ and $\bm B $
(analogues of block Jacobi methods)
or on their triangular decomposition (analogs of the Seidel block methods) 
\cite{VabishchevichAdditive}.
If necessary, apart from the main operator of the problem 
(the operator $\bm B$ in the equation (\ref{12}))
and the operator at the time derivative of the solution (operator $ \bm C $),
we will focus on the selection of the diagonal parts of the operator matrices $ \bm C $ and $ \bm B $ in splitting.
The splitting schemes themselves will be built by analogy with the work 
\cite{Vabishchevich2020}.

The minimal assumptions about the spaces $ V_i, \ i = 1, \dots, p $ are formulated as follows:
\begin{equation}\label{22}
 R_i R_i^* > 0, \quad  i = 1,\dots, p .
\end{equation} 
Select the diagonal parts of 
$\bm C$ and $\bm B$:
\[
 \bm C_0 = \mathrm{diag} \{R_1 R_1^*, \ldots, R_p R_p^* \} ,
 \quad \bm B_0 = \mathrm{diag} \{R_1 A R_1^*, \ldots, R_p A R_p^* \} .
\]
With the conditions
(\ref{4}) and (\ref{22}), we have
\begin{equation}\label{23}
  \bm C_0 = \bm C_0^* > 0,
  \quad \bm B_0 = \bm B_0^* > 0 .
\end{equation}

For the approximation
(\ref{12}), (\ref{13}), we will use three-level scheme
\begin{equation}\label{24}
\begin{split}
 \bm C_0 \left (\mu  \frac{\bm w^{n+1} - \bm w^{n}}{\tau } + (1-\mu) \frac{\bm w^{n} - \bm w^{n-1}}{\tau } \right ) & + 
 (\bm C - \bm C_0) \frac{\bm w^{n} - \bm w^{n-1}}{\tau } \\
 +  \bm B_0 (\sigma \bm w^{n+1} + (1-2\sigma)\bm w^{n} + \sigma\bm w^{n-1}) & +  (\bm B - \bm B_0) \bm w^{n}= \bm f^{n} , \\
 \quad n & = 1,\ldots,N-1, 
\end{split}
\end{equation} 
\begin{equation}\label{25}
 \bm w^{0} = \bm v^0 ,
 \quad \bm w^{1} = \widetilde{\bm v}^1 , 
\end{equation} 
with some weight parameters $ \mu = \mathrm {const}> 0$, 
$\sigma = \mathrm {const}> 0$.
The second initial condition (\ref{25}) is calculated, for example, 
using a two-level scheme (\ref{18}), (\ref{19}).
To find an approximate solution on a new level, 
we have a system of equations
\[
 (\mu \bm C_0  + \sigma \bm B_0) \bm w^{n+1} = \bm \varphi^n ,
\]  
which (see (\ref{23})) is uniquely solvable.
Stability conditions are given by the following statement.

\begin{thm}\label{t-1}

The three-level explicit-implicit scheme (\ref{22}) - (\ref{25}) 
is unconditionally stable for
\begin{equation}\label{26}
 \mu \geq  \frac{p}{2} ,
 \quad \sigma \geq \frac{p}{4} . 
\end{equation} 
Under these constraints, for an approximate solution to the problem (\ref{1}), 
(\ref{2}), the a priori estimate holds
\begin{equation}\label{27}
 \|y^{n+1/2}\|_A^2 \leq \|y^{1/2}\|_A^2  + 
 \left ( \bm D  \frac{\widetilde{\bm v}^1 - \bm v^0}{\tau }, \frac{\widetilde{\bm v}^1 - \bm v^0}{\tau } \right )
 + \frac{1}{2}  \sum_{k=1}^{n}\tau  \|f^k\|^2 ,
\end{equation} 
where the approximate solution to problem (\ref{1}), (\ref{2}) is
\[
 y^{n+1/2} = \frac{y^{n+1} + y^{n}}{2} ,
\] 
and the operator $\bm D = \bm D^* \geq 0$ has the representation
\begin{equation}\label{28}
 \bm D = \tau \left (\mu \bm C_0 - \frac{1}{2} \bm C \right )
 + \tau^2 \left (\sigma \bm B_0 - \frac{1}{4} \bm B \right ) . 
\end{equation} 
\end{thm}
\begin{pf}
Taking into account the equalities
\[
 \frac{\bm w^{n+1} - \bm w^{n}}{\tau } = \frac{\bm w^{n+1} - \bm w^{n-1}}{2\tau } +
 \frac{\tau}{2}  \frac{\bm w^{n+1} - 2 \bm w^{n} + \bm w^{n-1}}{\tau^2} ,
\] 
\[
 \frac{\bm w^{n} - \bm w^{n-1}}{\tau } = \frac{\bm w^{n+1} - \bm w^{n-1}}{2\tau } -
 \frac{\tau}{2}  \frac{\bm w^{n+1} - 2 \bm w^{n} + \bm w^{n-1}}{\tau^2} ,
\] 
\[
 \bm w^{n} = \frac{1}{4} (\bm w^{n+1} + 2 \bm w^{n} + \bm w^{n-1}) -
  \frac{\tau^2}{4}  \frac{\bm w^{n+1} - 2 \bm w^{n} + \bm w^{n-1}}{\tau^2} ,
\] 
the scheme (\ref{24}) can be written as
\begin{equation}\label{29}
\begin{split}
 \bm C \frac{\bm w^{n+1} - \bm w^{n-1}}{2\tau } & +
 \bm D \frac{\bm w^{n+1} - 2 \bm w^{n} + \bm w^{n-1}}{\tau^2}  \\
 & +
 \bm B \frac{1}{4} (\bm w^{n+1} + 2 \bm w^{n} + \bm w^{n-1}) = \bm f^{n} .
\end{split}
\end{equation} 

Taking into account
(\ref{14}) and inequality
\[
 \left ( \sum_{j=1}^{p} a_i \right )^2  \leq p \sum_{j=1}^{p} a^2_i ,
\] 
we have
\[
 (\bm C_0 \bm v, \bm v) =  \sum_{j=1}^{p} \left ( (R_i^* v_i)^2, 1 \right ) \geq   \frac{1}{p} (\bm C \bm v, \bm v) .
\] 
A similar relation holds for the operator $\bm B$ and its diagonal part 
$\bm B_0$, so that
\begin{equation}\label{30}
 \bm C \leq p \bm C_0,
 \quad \bm B \leq p \bm B_0 .
\end{equation} 
The operator $\bm D$, defined according to (\ref{29}) with (\ref{30}) 
and (\ref{26}), is non-negative: $\bm D = \bm D^* \geq 0 $.

We introduce the following new variables
\[
 \bm s^n = \frac{1}{2}(\bm w^{n} + \bm w^{n-1}) ,
 \quad \bm r^n = \bm w^{n} - \bm w^{n-1} .
\] 
The scheme (\ref{29}) can be re-written
\[
 \bm C \frac{\bm r^{n+1} + \bm r^{n}}{2\tau } +
 \frac{1}{\tau^2} \bm D(\bm r^{n+1} - \bm r^{n}) +
 \bm B \frac{\bm s^{n+1} + \bm s^{n}}{2}  = \bm f^{n} .
\]
We multiply this equation by (in a scalar way)
\[
 2(\bm s^{n+1} - \bm s^{n}) = \bm r^{n+1} + \bm r^{n-1} ,
\] 
and obtain the equality
\begin{equation}\label{31}
\begin{split}
 \frac{1}{2 \tau} (\bm C (\bm r^{n+1} + \bm r^{n})& ,\bm r^{n+1} + \bm r^{n})  +
 \frac{1}{\tau^2} (\bm D(\bm r^{n+1} - \bm r^{n}), \bm r^{n+1} + \bm r^{n}) \\ & +
 (\bm B (\bm s^{n+1} + \bm s^{n}), \bm s^{n+1} - \bm s^{n}) = (\bm f^{n}, \bm r^{n+1} + \bm r^{n})) . 
\end{split}
\end{equation} 
Considering (\ref{15}) and component-wise representation (see (\ref{10})) 
of the vector of the right-hand side
$\bm f^{n} = \{R_i f^n \} $, for two terms in (\ref{31}), we have
\[
 (\bm f^{n}, \bm r^{n+1} + \bm r^{n}) - \frac{1}{2 \tau} (\bm C (\bm r^{n+1} + \bm r^{n}) ,\bm r^{n+1} + \bm r^{n}) 
 \leq \frac{\tau }{2} \|f^n\|^2 .
\] 
This allow to obtain from 
(\ref{31}), the following inequality
\begin{equation}\label{32}
 \mathcal{E}^{n+1} \leq  \mathcal{E}^{n} + \frac{\tau }{2} \|f^n\|^2 ,
\end{equation} 
where
\[
 \mathcal{E}^{n+1} =  (\bm B \bm s^{n+1}, \bm s^{n+1}) + 
 \frac{1}{\tau^2} (\bm D(\bm r^{n+1}, \bm r^{n+1}) .
\] 
Taking into account that the operators
$\bm B$ and $\bm D$ are non-negative, and (\ref{32}), we get
\begin{equation}\label{33}
 \mathcal{E}^{n+1} \leq \mathcal{E}^{1} + \frac{1}{2}  \sum_{k=1}^{n}\tau  \|f^k\|^2 . 
\end{equation} 
Furthermore, taking into account 
(\ref{16}) and the introduced notations, we have 
\[
 \mathcal{E}^{n+1} \geq (\bm B \bm s^{n+1}, \bm s^{n+1}) = \left \| \frac{y^{n+1} + y^{n}}{2} \right \|_A^2 .
\] 
Thus, the estimate (\ref{27}) follows from (\ref{33}).
\end{pf}

According to Theorem \ref{t-1},
when using splitting scheme (\ref{22}) - (\ref{25}),
we find an approximate solution to the problem (\ref{1}) - (\ref{3})
at times $ t^{n + 1/2} $. But at the same time, we do not have 
a suitable stability estimate
for an approximate solution of the vector problem (\ref{12}), (\ref{13}).

\section{Solution decomposition on direct sum subspaces} 
\label{sec:4}
The main difficulties in constructing splitting schemes 
for a solution are due to
splitting the operator $\bm C$ for the time derivatives of 
the vector $\bm v $ in the equation (\ref{12}).
For this reason, we need to focus on three-level schemes.
We note if these difficulties can be avoided, one can 
 use two-level schemes. One can achieve a
success by narrowing the class of the family of
spaces $ V_i, \ i = 1, \ldots, p $, by a stricter 
choice of restriction operators $ R_i, \ i = 1, \ldots, p $.

Let the space 
$U$ be the direct sum of subspaces
$V_i, \ i = 1,\ldots, p$:
\[
 U = R_1^* V_1 \oplus \ldots \oplus R_p^* V_p,
\] 
such that for every $ u \in U $ there is a representation (\ref{9}).
In this case
\[
 v_i = R_i u, 
 \quad i = 1,\ldots, p , 
\] 
and therefore for the 
restriction operators, we have
\begin{equation}\label{34}
 R_i R_j^* = \left \{
 \begin{array}{cc}
  I_i ,  &  i=j, \\
  0 ,  &  i \neq j , \\
\end{array}
\right .  
\end{equation} 
where $I_i$ --- unit operator in $V_i, \ i = 1,\ldots, p$.
Previously, weaker constraints (\ref{22}) were used instead of (\ref{34})
for restriction operators $R_i, \ i = 1, \ldots, p$.

For (\ref{34}), the Cauchy problem is solved for the equation
\begin{equation}\label{35}
 \frac{d \bm v}{d t} + \bm B \bm v = \bm f ,
\end{equation} 
i.e., in (\ref{12}) $\bm C = \bm C_0 = \bm I$. 
We can restrict ourselves to heterogeneous
approximation of only the operator $\bm B$. 
Let us note  new possibilities of constructing 
splitting for an approximate solution of the problem (\ref{13}), (\ref{35}), 
connecting them, in particular, with the schemes
of second order time approximation.

By analogy with (\ref{24}), for an approximate solution of the problem 
(\ref{13}), (\ref{35}), we will use the scheme
separating the diagonal part of the operator $\bm B$:
\begin{equation}\label{36}
\begin{split}
 \frac{\bm w^{n+1} - \bm w^{n-1}}{2\tau } & +
 \bm B_0 (\sigma \bm w^{n+1} + (1-2\sigma)\bm w^{n} + \sigma\bm w^{n-1}) \\
 & +  (\bm B - \bm B_0) \bm w^{n}= \bm f^{n} , 
 \quad n = 1,\ldots,N-1 .
\end{split}
\end{equation} 

\begin{thm}\label{t-2}

The three-level scheme of the second order approximation (\ref{25}), 
(\ref{34}), (\ref{36}) is unconditionally stable for $\sigma \geq p/4$, 
and for the approximate solution of the problem (\ref{1}), (\ref{2}), 
 a priori estimate (\ref{27}) holds where
\[
 \bm D = \tau^2 \left (\sigma \bm B_0 - \frac{1}{4} \bm B \right ) . 
\]
\end{thm}
\begin{pf}
The proof is carried out 
similarly to the proof of Theorem \ref{t-1}.
\end{pf}

A simpler variant of the splitting scheme with the selection of the diagonal part of the problem operator is associated with the use of a two-level scheme
\begin{equation}\label{37}
\begin{split}
 \frac{\bm w^{n+1} - \bm w^{n}}{\tau } & +
 \bm B_0 (\sigma \bm w^{n+1} + (1-\sigma)\bm w^{n}) \\
 & +  (\bm B - \bm B_0) \bm w^{n}= \bm f^{n} , 
 \quad n = 0,\ldots,N-1 .
\end{split}
\end{equation} 
The stability conditions for this scheme are given by the following statement.

\begin{thm}\label{t-3}
The two-level scheme of the first order approximation (\ref{19}), (\ref{34}), 
(\ref{37}) is unconditionally stable for $\sigma \geq p/2 $, 
and for the approximate solution of the problem (\ref{1}), (\ref{2}),
  a priori estimate (\ref{7}) holds.
\end{thm}
\begin{pf}
We rewrite 
(\ref{37}) in the form
\begin{equation}\label{38}
 \left ( \bm I + \tau \left (\sigma \bm B_0 - \frac{1}{2} \bm B \right )  \right ) \frac{\bm w^{n+1} - \bm w^{n}}{\tau }  +
 \bm B \frac{\bm w^{n+1} + \bm w^{n}}{2}  = \bm f^{n} ,   
\end{equation} 
Taking into account the inequality
(\ref{30}) with $\sigma \geq p/2$, we have
\[
 \sigma \bm B_0 - \frac{1}{2} \bm B \geq 0 .
\] 
Multipying
(\ref{38}) (scalarly) by $2(\bm w^{n+1} - \bm w^{n})$, we have
\[
 2 \tau \left \| \frac{\bm w^{n+1} - \bm w^{n}}{\tau }\right \|^2 + (\bm B \bm w^{n+1}, \bm w^{n+1}) - (\bm B \bm w^{n}, \bm w^{n})
 \leq  2  (\bm f^{n}, \bm w^{n+1} - \bm w^{n}) .
\] 
For the estimate on the right hand side, we use 
\[
 (\bm f^{n}, \bm w^{n+1} - \bm w^{n}) \leq \tau \left \| \frac{\bm w^{n+1} - \bm w^{n}}{\tau }\right \|^2
 + \frac{\tau }{4} \|\bm f^{n}\|^2 ,
\] 
and therefore
\[
 (\bm B \bm w^{n+1}, \bm w^{n+1}) \leq  (\bm B \bm w^{n}, \bm w^{n}) + \frac{\tau }{2} \|\bm f^{n}\|^2 .
\] 
This inequality immediately implies the estimate (\ref{7}) 
for the approximate solution (\ref{1}), (\ref{2}).
\end{pf}

Separately, we single out the class of splitting schemes for the problem 
(\ref{13}), (\ref{35}), which is associated with triangular splitting
operator $\bm B$. We define
\[
 \bm B_1 = \left (\begin{array}{cccc}
  \frac{1}{2} R_1 R R_1^*  &  0 & \cdots &  0 \\
  R_2 R R_1^*   &  \frac{1}{2} R_2 R R_2^*  & \cdots &  0 \\
  \cdots  & \cdots & \cdots &  0 \\
   R_p R R_1^*  &  R_p R R_2^*  & \cdots &  \frac{1}{2} R_p R R_p^*  \\
\end{array}
 \right ) ,
\]
\[
  \bm B_2 = \left (\begin{array}{cccc}
  \frac{1}{2} R_1 A R_1^*  &  R_1 A R_2^*  & \cdots &  R_1 A R_p^*  \\
  0  &  \frac{1}{2} R_2 A R_2^*  & \cdots &  R_2 A R_p^*  \\
  0  & \cdots & \cdots & \cdots \\
  0  &  0  & \cdots &  \frac{1}{2} R_p R R_p^*  \\
\end{array}
 \right ) ,
\] 
and, thus,
\begin{equation}\label{39}
 \bm B = \bm B_1 +  \bm B_2 ,
 \quad  \bm B_1^*  =  \bm B_2 .
\end{equation}

For two-component splitting (\ref{39}), more interesting cases are
factorized schemes
\cite{SamarskiiTheory, VabishchevichAdditive}.
For an approximate solution of the problem (\ref{13}), (\ref{35}), (\ref{39}),
we will use the scheme
\begin{equation}\label{40}
\begin{split}
 (\bm I + \tau \sigma \bm B_1) (\bm I + \tau \sigma \bm B_2) \frac{\bm w^{n+1} - \bm w^{n}}{\tau }
 +  \bm B\bm w^{n}= \bm f^{n} , 
 \quad n = 0,\ldots,N-1 .
\end{split}
\end{equation}    
For $\sigma = 1$, the factorized scheme (\ref{40}) 
has the first order  approximation in time and
has a similar operator analog 
as the classical Douglas-Rachford scheme.
Of particular interest is the factorized scheme (\ref{40}) 
of the second order approximation,
when $\sigma = 1/2 $, it is 
an analogue of the Peaceman-Rachford scheme.

\begin{thm}\label{t-4}
The factorized scheme (\ref{19}), (\ref{39}), (\ref{40}) is unconditionally 
stable for $\sigma \geq 1/2 $, and for an approximate solution of 
the problem (\ref{1}) , (\ref{2}), a priori estimate (\ref{7}) holds.
\end{thm}
\begin{pf}
Similar to 
(\ref{38}), we rewrite (\ref{40}) in the form
\[
 (\bm I + \tau \bm G)  \frac{\bm w^{n+1} - \bm w^{n}}{\tau }  +
 \bm B \frac{\bm w^{n+1} + \bm w^{n}}{2}  = \bm f^{n} ,    
\] 
where
\[
 \bm G = \left (\sigma - \frac{1}{2} \right ) \bm B + \sigma^2 \tau \bm B_1 \bm B_2 .   
\] 
Taking into account  (\ref{39}) with $\sigma \geq 1/2$, we have 
\[
 \bm G = \bm G^* \geq 0 .
\] 
The rest of the proof is similar to the proof of Theorem 
\ref{t-3}.

\end{pf}

\section{Conclusions}

In this paper, we discuss splitting methods and solution decomposition
and investigate splitting methods for solution decomposition.
The proposed methods, in general, can 
combine spatial decomposition of the solution and temporal splitting 
for efficient solution computations of nonstationary problems.  
The proposed splitting method results to equations for
 individual solution components that
 consist of the time derivatives of the solution components.
Unconditionally stable splitting schemes are constructed and analyzed
for a first-order evolution equation formulated
 in a finite-dimensional Hilbert space. 
In our splitting algorithms, we consider  the decomposition
 of both the main operator
of the system and the operator at the time derivative.
These approaches are motivated by various applications, 
including model reduction, multiscale methods, 
domain decomposition, 
and so on, where spatially local problems are 
solved and used to compute the solution. This paper presents a general framework that can be used to guide, for example,
 spatial decomposition approaches in non-stationary problems.

\section*{Acknowledgements}

The publication has been prepared with the support by the mega-grant of the Russian Federation Government 14.Y26.31.0013.

\end{document}